\documentclass[preprint,3p,natbib]{elsarticle}

\biboptions{sort&compress}
\usepackage{amsmath}
\usepackage{mathrsfs}
\usepackage{amssymb}
\usepackage{amsthm}
\usepackage{makecell}
\usepackage{multirow}
\usepackage{makecell}
\usepackage{diagbox}
\usepackage{graphicx}
\usepackage{float}
\usepackage{xcolor}
\usepackage{xparse}
\usepackage{multirow}
\ExplSyntaxOn
\DeclareExpandableDocumentCommand{\IfNoValueOrEmptyTF}{mmm}
 {
  \IfNoValueTF{#1}{#2}
   {
    \tl_if_empty:nTF {#1} {#2} {#3}
   }
 }
\ExplSyntaxOff
\newcommand{\mathttul}[3][\mkern2mu]{\underline{\mathtt{#2}_{#3}\mkern-2mu}^{#1}}
\DeclareDocumentCommand\mathttul{o o m m}{
\IfNoValueOrEmptyTF{#2}{
\underline{\mathtt{#3}_{#4}\mkern-2mu}\IfNoValueF{#1}{^{#1}}\mkern3mu
}{
(\underline{\mathtt{#3}_{#4}\mkern-2mu})^{#1}_{#2}
}
}

\newtheorem{example}{Example}[section]
\newtheorem{lemma}{Lemma}[section]
\newtheorem{theorem}{Theorem}[section]

\newtheorem{remark}{Remark}[section]

\numberwithin{equation}{section}

\journal{JMAA}

\begin{document}

\begin{frontmatter}



\title{\textbf{Quantitative recurrence problem on some Bedford-McMullen carpets}}


\author[1]{Yu-Liang Wu}
\ead{yu-liang.wu@oulu.fi}
\author[2]{Na Yuan\corref{*}}
\ead{yuanna199@gmail.com}

\affiliation[1]{organization={Department of Mathematical Sciences, University of Oulu},
addressline={Pentti Kaiteran katu 1},
            city={Oulu},
            postcode={90014}, 
            country={Finland}}
 \affiliation[2]{organization={School of Mathematics, Guangdong University of Education},
 addressline={Xingang Middle Road 351},
            city={Guangzhou},
            postcode={510303}, 
            country={China}}
\cortext[*]{Corresponding author}

\begin{abstract}
In this paper, we study the Hausdorff dimension of the quantitative recurrent set of the canonical endomorphism on the Bedford–McMullen carpets whose Hausdorff dimension and box dimension are equal.
\end{abstract}


\end{frontmatter}


\section{Introduction}

\subsection{Background}
The concept of recurrence plays an important role in dynamical systems and ergodic theory. Let $(X,\mathcal{B},\mu,T)$ be a measure-preserving system equipped with a compatible metric $d$, i.e., $(X,d)$ is a metric space and $\mathcal{B}$ is a Borel $\sigma$-algebra of $X$. If $(X, d)$ is a separable metric space, the Poincar\'{e} Recurrence Theorem implies that $\mu$-almost every $x\in X$ is recurrent in the sense that
\begin{equation*}
    \liminf_{n\rightarrow\infty}d(T^nx, x)=0.
\end{equation*}
In nature, the result provides no information about the rate at which an orbit will return to the initial point or in what manner a neighborhood of the initial point will shrink under the iteration. The interest in these quantitative characterization has provoked a rich subsequent literature on the so-called quantitative recurrent sets: given a \emph{rate function} $\psi:X\times \mathbb{N} \rightarrow (0,\infty)$, the \emph{quantitative recurrent set with respect to $\psi$} is defined as
\begin{equation} \label{eq:quantitative_recurrent_set}
    R(T,\psi)=\big\{x\in X:d(T^nx, x)<\psi(n,x) \text{ for infinitely many } n \in \mathbb{N} \big\}.
\end{equation}
Boshernitzan \cite{Boshernitzan1993} gave an outstanding result for general systems concerning the size in measure of $R(T,\psi)$. Later, Barreira and Saussol \cite{Barreira2001} stated a finer result. 

In recent years, many authors have turned their eyes to recurrent sets on fractals. On the one hand, some researchers showed that the $\mu$-measure of the set $R(T,\psi)$ is null or full according to convergence or divergence of a certain series in some dynamical systems (see Chang-Wu-Wu \cite{Chang2018}, Baker-Farmer \cite{Baker2019}, Hussain-Li-Simmons-Wang \cite{Hussain2021}, Kirsebom-Kunde-Persson \cite{Kirsebom2022}, Persson \cite{Persson2022}, Kleinbock-Zheng \cite{Kleinbock2023} and Baker-Koivusalo \cite{Baker2024}). On the other hand, many researchers studied the Hausdorff dimension of the set $R(T,\psi)$ in some dynamical systems (see Tan-Wang \cite{Tan2011} and Seuret-Wang \cite{Seuret2013}). Note that when we require $\{T^n x\}_{n\ge 1}$ to return to the neighborhoods of a chosen point $x_0\in X$ rather than the initial point $x$, the problem becomes the so-called shrinking target problem, which was first investigated by Hill and Velani \cite{Hill1995}. Since then, many authors have contributed to the study of the shrinking target problem. To name but a few, see \cite{Allen2021ONTH, Barany2018ShrinkingTO, Bugeaud2003, Fan2013AMM, Hill2002, koivusalo2018, li2022, Li2014, liao_seuret_2013, persson_rams_2017, Shen2013ShrinkingTP, Tseng2007OnCR} and references within.

As for the Hausdorff dimension, it is to be noted that the aforementioned works mainly involve systems of $\mathbb{R}^1$ or dynamical systems in conformal dynamics, and hardly anything is known as far as high-dimensional non-conformal dynamics are concerned. The only known result was presented by B\'{a}r\'{a}ny and Troscheit \cite{Barany2021}, who investigated dimensions for both the quantitative recurrence and the shrinking target problems for dynamically defined subsets of generic self-affine sets (in the sense of Lebesgue almost all translations). In this work, we intend to consider the Bedford-McMullen carpets, a typical family of planar self-affine sets introduced in \cite{bedfordcrinkly} and \cite{McMullen1984}, and discuss the recurrence based on the Euclidean metric. This is the first time the problem of recurrent sets for deterministic self-affine fractal sets is discussed. We would like to point out that due to the similarity of the quantitative recurrence problem and the shrinking target problem in nature, techniques for one problem sometimes apply to the other. For our study, arguments closely follow those from Barany and Rams \cite{Barany2018ShrinkingTO} with a technical nuance mentioned in Remark \ref{rem:comparison_barany}.

\subsection{The main theorem}
It is the purpose of this paper to study the quantitative recurrent set problem on the system $(K, T)$ of Bedford-McMullen carpet $K$. 
Let $2 \le m_1 \le m_2$ be two integers and $\Sigma_{m_i}=\{0,1,\cdots,m_i-1\}$ for $i=1, 2$. Define for every $\mathtt{a}=(\mathtt{a}^{(1)},\mathtt{a}^{(2)}) \in \Sigma_{m_1} \times \Sigma_{m_2}$ a map $\phi_{\mathtt{a}}:[0,1]^2 \to [0,1]^2$ as
\[
    \phi_{\mathtt{a}}(x^{(1)},x^{(2)}):=\left(\frac{x^{(1)}+\mathtt{a}^{(1)}}{m_1},\frac{x^{(2)}+\mathtt{a}^{(2)}}{m_2}\right).
\]
Given any nonempty subset $A \subseteq \Sigma_{m_1} \times \Sigma_{m_2}$, the \emph{Bedford-McMullen carpet $K$ associated with $A$} is the unique attractor of the iterated function system $\{\phi_{\mathtt{a}}: \mathtt{a} \in A\}$. If we consider the coding map $\pi: (\Sigma_{m_1} \times \Sigma_{m_2})^{\mathbb{N}} \to [0,1]^2$ defined as
\[
    \pi(\mathtt{x})=\left(\sum_{n=1}^{\infty}\frac{\mathtt{x}^{(1)}_n}{m_1^n},\sum_{n=1}^{\infty}\frac{\mathtt{x}^{(2)}_n}{m_2^n}\right),
\]
where $\mathtt{x}=(\mathtt{x}^{(1)}_1,\mathtt{x}^{(2)}_1)(\mathtt{x}^{(1)}_2,\mathtt{x}^{(2)}_2)(\mathtt{x}^{(1)}_3,\mathtt{x}^{(2)}_3) \cdots$,
then the Bedford-McMullen carpet $K$ can be expressed as $K=\pi(A^{\mathbb{N}})$. It is noteworthy that the latter definition naturally endows the carpet $K$ with a map $T:K \to K$ defined as
\[
    T(x)=(T_{m_1}(x^{(1)}), T_{m_2}(x^{(2)})):=\left(m_1 x^{(1)}\ (\bmod\ 1),m_2 x^{(2)}\ (\bmod\ 1)\right).
\]

Let $\psi: \mathbb{N} \to (0,\infty)$ be a \emph{rate function} and $\gamma > 0$, and define the \emph{recurrent set $W_{\gamma}(K,T,\psi)$ with respect to $\psi$} as:
\begin{align}\label{e2}
    W_{\gamma}(K,T,\psi):=\left\{x \in K: \begin{cases}
        |x^{(1)}-T^n_{m_1}(x^{(1)})| < \psi(n) \\
        |x^{(2)}-T^n_{m_2}(x^{(2)})| < \psi(n)^{\gamma} 
    \end{cases} \text{for infinitely many } n \right\}.
\end{align}
Note that for the rest of the discussion, we assume no monotonicity of $\psi$. Denote by $\dim_{\mathrm{H}}$ and $\dim_{\mathrm{B}}$ the Hausdorff and box dimensions, respectively. The following notations are used throughout the discussions.
\begin{equation}\label{E0}
    \ell_{1}(n) = -\log_{m_1} \psi(n) \quad \text{and} \quad \ell_{2}(n) = -\log_{m_2} \psi(n)^{\gamma}
\end{equation}
and
\begin{equation}\label{E1}
    \frac{\ell_{i}(n)}{n} =\tau_i(n) \quad \text{and} \quad \liminf_{n \to \infty}\frac{\ell_{i}(n)}{n}=\tau_i, \ \text{ for } i=1,2.
\end{equation}
Denote
\begin{align}\label{E8}
    \hat{\ell}_{i}(n)=\lceil \ell_{i}(n) \rceil := \min\{k\in \mathbb{N}: k\geq \ell_{i}(n)\}.
\end{align}
for $i=1,2$ and $n\in \mathbb{N}$. It is noteworthy that if $0 \leq \gamma \leq \log_{m_1}m_2$, then $\ell_1(n) \ge \ell_2(n)$ and $\hat{\ell}_1(n) \ge \hat{\ell}_2(n)$ for all $n \in \mathbb{N}$.

Our main result related to the Hausdorff dimension of the set $W_{\gamma}(K,T,\psi)$ is as follows.
\begin{theorem}\label{M1}\label{M1-1}
    Let $ W_{\gamma}(K,T,\psi)$ be defined in \eqref{e2}, $\tau_1 \ge 0$, and $0\leq \gamma\leq \log_{m_1}m_2$. Let $M$ be the number of columns containing at least one chosen rectangle and $N_i$ the number of rectangles chosen from the $i$-th non-empty column.
    Suppose that $\dim_{\mathrm{B}} K = \dim_{\mathrm{H}} K$, or equivalently, there exists an integer $N\geq 1$ such that $N_i \in \{0, N\}$ for all $i$ (see \cite{bedfordcrinkly} or \cite{McMullen1984}), we have $\dim_{\mathrm{H}} W_{\gamma}(K,T,\psi) = \min \{t_1,t_2\}$, where
    \[
        t_1= \begin{cases}
            \left(1-\frac{\tau_1 \log_{m_2} m_1}{1+\tau_2}\right) \cdot \log_{m_1} M + \frac{\log_{m_2} N}{1+\tau_2} & \text{if } 1+(1-\gamma)\tau_1 \le \log_{m_1} m_2, \\
            \frac{\log_{m_2} M + \log_{m_2} N}{1+\tau_2} & \text{if } 1+(1-\gamma)\tau_1 > \log_{m_1} m_2.
        \end{cases}
    \]
    \[
        t_2 = \begin{cases}
            \frac{\log_{m_1} M}{1+\tau_1} + \log_{m_2} N & \text{if } 1 + \tau_1 \le \log_{m_1} m_2, \\
            \frac{\log_{m_1} M + \log_{m_1} N}{1+\tau_1} & \text{if } 1 + (1 - \gamma) \tau_1 \le \log_{m_1} m_2 < 1 + \tau_1, \\
            \frac{\log_{m_1} M + (1 + (1 - \gamma) \tau_1) \log_{m_2} N}{1+\tau_1} & \text{if } 1+(1-\gamma) \tau_1 > \log_{m_1} m_2.
        \end{cases}
    \]
    In particular, if $\gamma = 1$, this is reduced to the following form.
    
    $(1)$ If $\log_{m_1} m_2 > 1 + \tau_1$, then
    \begin{align*}
    \dim_{\mathrm{H}} W_{\gamma}(K,T,\psi)=
	\min \left\{\frac{\log_{m_1} M +  \log_{m_2} N}{1+\tau_2},
        \frac{\log_{m_1} M}{1+\tau_1}  + \log_{m_2} N\right\}.
    \end{align*}
    
    $(2)$ If  $\log_{m_1} m_2 \le 1 + \tau_1$, then 
    \begin{align*}
        \dim_{\mathrm{H}} W_{\gamma}(K,T,\psi)  = \min\left\{
        \frac{\log_{m_1} M +\log_{m_2}N}{1+\tau_2}, \
          \frac{\log_{m_1} M +\log_{m_1} N}{1+\tau_1}   \right\},
    \end{align*}
\end{theorem}
\begin{remark}\label{remark-1}
    We note that Theorem \ref{M1-1} covers the cases of all nontrivial $\tau_1$ for if $\tau_1<0$, then $W_{\gamma}(K,T,\psi) = K$. 
    In fact, for the rest of the discussion, we further assume $0 < \tau_1 < \infty$ since the cases $\tau_1 = 0$ and $\tau_1 = \infty$ could be inferred from the remaining cases by noting that $W_{\gamma}(K,T,\psi) \subset W_{\gamma}(K,T,\phi)$ if $\psi \le \phi$.
\end{remark}
\begin{remark}
    The set $W_1(K,T,\psi)$ can be interpreted as a quantitative recurrent set defined in \eqref{eq:quantitative_recurrent_set} with the maximum norm of $\mathbb{R}^2$. Nevertheless, Theorem \ref{M1} implies that the dimension of the quantitative recurrent set is invariant under all equivalent metrics, including the one induced by the Euclidean norm.
\end{remark}
\begin{remark}
    The homogeneity property (i.e., $N_i \in \{0, N\}$ for all $i$) in Theorem \ref{M1} is imposed only for the lower bound (derived in Section 3) and upper bound (derived in Section 4) to coincide. Nevertheless, our arguments for the lower and upper bounds remain valid (with necessary modifications) without this assumption. As a related issue, in Theorem \ref{M1-1}, our study of $W_{\gamma}(K,T,\psi)$ is limited to a restrictive class of $\gamma$ so that $\tau_1 \ge \tau_2$ and that the homogeneity property could be exploited to deduce matching lower bound and upper bound. We note that this is also the main difficulty we face if we were to further generalize the theorem by having independent rate functions in the two coordinates. 
\end{remark}

\subsection{Organization of the paper}
In this paper, we will present a Hausdorff dimension formula to the recurrent set valid for the subfamily of Bedford–McMullen carpets whose Hausdorff dimensions are equal to the box dimensions. The paper is organized as follows. In Section \ref{S1}, we present necessary notation and preliminaries. Section \ref{S2} and Section \ref{S3} are devoted to the proof of Theorem \ref{M1}. Finally, in Section \ref{S4}, we provide some examples to which the main theorem is applicable.
\section{Preliminaries}\label{S1}
In this section, we introduce necessary notations and preliminaries. For convenience's sake, the following notations of the symbolic space are also introduced. For $i=1,2$, let $\Sigma_{m_i}^n=\{\mathtt{u}:\mathtt{u}=(\mathtt{u}_1,\dots,\mathtt{u}_n), \mathtt{u}_j\in \Sigma_{m_i}, \ j=1,\dots,n \}$. Denote by $|u|$ the length of $u\in \Sigma_{m_i}^n$ for $i=1,2$. Firstly, for any $l\in \mathbb{N} \cup \{\infty\}$ and $u \in \Sigma_{m_i}^n$, write $(u)^l$ for the word $(u,\cdots,u)$ ($l$ times repeated concatenation of the word). More generally, for any positive number $l>0$, denote by $(u)^{l}$ the word $(u)^{ \lfloor l \rfloor }u'$, where $u'$ is the prefix of $u$ with length $\lfloor(l-\lfloor l \rfloor)|u|\rfloor$ and $\lfloor l \rfloor=\max\{k\in \mathbb{N}: k\leq l\}$. Secondly, we identify the spaces $(\Sigma_{m_1} \times \Sigma_{m_2})^{\mathbb{N}}$ and $(\Sigma_{m_1}^{\mathbb{N}} \times \Sigma_{m_2}^{\mathbb{N}})$ by setting 
\[
    (\mathtt{x}^{(1)}_1,\mathtt{x}^{(2)}_1)(\mathtt{x}^{(1)}_2,\mathtt{x}^{(2)}_2)(\mathtt{x}^{(1)}_3,\mathtt{x}^{(2)}_3) \cdots \sim (\mathtt{x}^{(1)}_1 \mathtt{x}^{(1)}_2 \mathtt{x}^{(1)}_3 \cdots ,\mathtt{x}^{(2)}_1 \mathtt{x}^{(2)}_2 \mathtt{x}^{(2)}_3 \cdots),
\]
which is a one-to-one correspondence. With this identification, we introduce the following notation of cylinder set with different lengths in its coordinates:
\[
    [(\mathtt{w}^{(1)},\mathtt{w}^{(2)})]=\big\{(\mathtt{x}^{(1)},\mathtt{x}^{(2)}): \mathtt{x}^{(i)}_{1:n_i} = \mathtt{w}^{(i)}, i=1,2\big\}, \quad \big(\mathtt{w}^{(1)},\mathtt{w}^{(2)}\big) \in \Sigma_{m_1}^{n_1} \times \Sigma_{m_2}^{n_2}\big\},
\]
each of which is associated with a half-open half-closed rectangle
\[
    I([(\mathtt{w}^{(1)},\mathtt{w}^{(2)})])=\left[\sum_{i=1}^{n_1} \frac{\mathtt{w}^{(1)}}{m_{1}^i},\sum_{i=1}^{n_1} \frac{\mathtt{w}^{(1)}}{m_{1}^i}+m_{1}^{-n_1}\right) \times \left[\sum_{i=1}^{n_2} \frac{\mathtt{w}^{(2)}}{m_{2}^i},\sum_{i=1}^{n_2} \frac{\mathtt{w}^{(2)}}{m_{2}^i}+m_{2}^{-n_2}\right).
\]
We call $I([(\mathtt{w}^{(1)},\mathtt{w}^{(2)})])$ an \emph{$n_1$-th level rectangle} when $n_1=n_2$. In particular, we call $I([(\mathtt{w}^{(1)},\mathtt{w}^{(2)})])$ an \emph{$n_1$-th level approximate square} when $n_2 = \lceil \log_{m_2} m_1 \cdot n_1 \rceil$. For any $\mathtt{x}\in (\Sigma_{m_1} \times \Sigma_{m_2})^{\mathbb{N}}$, let $\mathtt{x}^{(i)}_{1:n}:=\mathtt{x}^{(i)}_{1}\mathtt{x}^{(i)}_{2}\dots\mathtt{x}^{(i)}_{n}$. If $A\subset \Sigma_{m_1} \times \Sigma_{m_2}$, let 
\[
    A^n=\big\{(\mathtt{x}^{(1)}_1 \mathtt{x}^{(1)}_2 \mathtt{x}^{(1)}_3 \cdots \mathtt{x}^{(1)}_n,\mathtt{x}^{(2)}_1 \mathtt{x}^{(2)}_2 \mathtt{x}^{(2)}_3 \cdots \mathtt{x}^{(2)}_n): (\mathtt{x}^{(1)}_i,\mathtt{x}^{(2)}_i)\in A\big\}.
\]
For any $\mathtt{u}=(\mathtt{u}_1^{(1)}\mathtt{u}_2^{(1)}\dots\mathtt{u}_n^{(1)}, \mathtt{u}_1^{(2)}\mathtt{u}_2^{(2)}\dots\mathtt{u}_n^{(2)})\in A^n$, let $\mathtt{u}^{(i)}=(\mathtt{u}_1^{(i)}\mathtt{u}_2^{(i)}\dots\mathtt{u}_n^{(i)} )$, $i=1,2$.

\section{Proof of the lower bound}\label{S2}

In this section, we will prove the lower bound of the Hausdorff dimension of $W_{\gamma}(K,T,\psi)$, and to this end, we are going to prove Lemma \ref{lem:Markov}-\ref{L5}.

In the following, we outline our strategy of proof. To begin with, we note that it suffices to prove the lower bound of Theorem \ref{M1-1} with an additional assumption that $0 \le \hat{\ell}_2(n) \le \hat{\ell}_1(n)$ for all $n$. Indeed, if Theorem \ref{M1-1} holds under the assumption, then for any $\psi(n)$ one may choose $\phi(n) = \min\{1, \psi(n)\}$ so that $\dim_{\mathrm{H}} W_{\gamma}(K,T,\psi) \ge \dim_{\mathrm{H}} W_{\gamma}(K,T,\phi)$ yields the desired lower bound. Let $p \in [0,1]^{A}$ be a probability vector indexed by $A$. For every $\mathtt{a} \in A$, let 
\[
    p_{\mathtt{a}^{(1)}}=
    \sum_{\mathtt{b}\in A, \ \mathtt{b}^{(1)}=\mathtt{a}^{(1)}}p_{\mathtt{b}} \hspace{2em} \text{and} \hspace{2em}
    A^{(1)}=\{\mathtt{a}^{(1)}: \mathtt{a}=(\mathtt{a}^{(1)}, \mathtt{a}^{(2)}), \ \mathtt{a} \in A\}. 
\]
We take an increasing sequence of natural numbers $n_i$ such that 
\begin{align}\label{E30}
    2^i \sum_{j=1}^{i}n_j\ll n_{i+1} \quad \text{ and } \quad \lim_{i \to \infty} \tau_1(n_j) = \tau_1,
\end{align}
and define the measure $\mu$ (with respect to probability measure $p$) as follows. To begin with, write $\mathcal{W}_1 = A^{n_1}$ and set
\begin{align*}
    \mathcal{W}_{i+1} = \left\{\mathtt{w} \mathtt{u} \mathtt{v} \mathtt{w}'\in A^{n_{i+1}}: \mathtt{w} \in \mathcal{W}_{i}, \mathtt{u} \in A^{\hat{\ell}_2(n_i)}, \mathtt{v} \in A^{\hat{\ell}_1(n_i)-\hat{\ell}_2(n_i)}, \mathtt{w}' \in A^{n_{i+1}-(n_i+\hat{\ell}_1(n_i))}, \right. \\
    \left. (\mathtt{u} \mathtt{v})^{(1)} = (\mathtt{w}^{(1)})^{\frac{\hat{\ell}_1(n_{i})}{n_{i}}}, \mathtt{u}^{(2)} = (\mathtt{w}^{(2)})^{\frac{\hat{\ell}_2(n_{i})}{n_{i}}}\right\}
\end{align*}
Observe that by construction, a word $\mathtt{w}$ is in $\mathcal{W}_i$ if and only if it can be written in the form of
\begin{equation} \label{eq:decomposition}
    \mathtt{w} = \mathttul{w}{1} \mathttul{u}{1} \mathttul{v}{1} \mathttul{w}{2} \cdots \mathttul{w}{i-1} \mathttul{u}{i-1} \mathttul{v}{i-1} \mathttul{w}{i},
\end{equation}
where, by setting $n_0 = 0$, $\mathttul{w}{j} \in A^{n_{j}-(n_{j-1}+\hat{\ell}_1(n_{j-1}))}$ is arbitrary, and $\mathttul{u}{j} \in A^{\hat{\ell}_2(n_{j})}$, $\mathttul{v}{j} \in A^{\hat{\ell}_1(n_{j-1}) - \hat{\ell}_2(n_{j-1})}$ satisfy
\begin{equation} \label{eq:recurrence_assumption}
\begin{cases}
    \mathttul[(1)]{u}{j} \mathttul[(1)]{v}{j} = \left(\mathttul[(1)]{w}{1} \cdots \mathttul[(1)]{u}{j-1} \mathttul[(1)]{v}{j-1} \mathttul[(1)]{w}{j}\right)^{\frac{\hat{\ell}_{1}(n_j)}{n_j}}, \\
    \mathttul[(2)]{u}{j} = \left(\mathttul[(2)]{w}{1} \cdots \mathttul[(2)]{u}{j-1} \mathttul[(2)]{v}{j-1} \mathttul[(2)]{w}{j}\right)^{\frac{\hat{\ell}_{2}(n_j)}{n_j}},
\end{cases}
\end{equation}
for all $1 \le j \le i$. For convenience, denote
\begin{align*}
    L_{i} &= |\mathttul{w}{1} \cdots \mathttul{u}{i}| = n_i + \hat{\ell}_2(n_i), \\
    R_{i} &= |\mathttul{w}{1} \cdots \mathttul{v}{i}| = n_i + \hat{\ell}_1(n_i).
\end{align*}
Next, we inductively define $\mu$ on cylinder sets of length $n_i+\hat{\ell}_1(n_i)$ for $i \in \mathbb{N}$ as 
\begin{equation} \label{eq:prob_measure}
    \mu[\mathtt{w}] = \begin{cases}
        \mu[\mathtt{w}'] \prod\limits_{\ell=1}^{|\mathtt{v}|} \frac{p_{\mathtt{v}_{\ell}}}{p_{\mathtt{v}^{(1)}_{\ell}}} \prod\limits_{\ell=1}^{|\mathtt{w}''|} p_{\mathtt{w}''_{\ell}};
        & \text{if } \mathtt{w}' \in \mathcal{W}_{i-1} \text{ and } \mathtt{w} = \mathtt{w}' \mathtt{u} \mathtt{v} \mathtt{w}'' \in \mathcal{W}_i \\
        0 & \text{otherwise.}
    \end{cases}
\end{equation}
It is readily checked by definition that 
\[
    1 = \sum_{\mathtt{w}' \in A^{n_{1}+\hat{\ell}_1(n_{1})}} \mu[\mathtt{w}'] \quad \text{and} \quad \mu[\mathtt{w}] = \sum_{\mathtt{w}' \in A^{n_{i+1}+\hat{\ell}_1(n_{i+1})}} \mu[\mathtt{w}'] \text{ for every } \mathtt{w} \in A^{n_i+\hat{\ell}_1(n_i)},
\]
which extends $\mu$ to an additive set function on the algebra generated by all cylinders sets and hence, by Carath\'{e}odory extension theorem, to a Borel probability measure on $A^{\mathbb{N}}$. It then follows from this definition that $\text{supp} (\mu \circ \pi^{-1})\subseteq W_{\gamma}(K,T,\psi)$. We note that our choice of measure is similar to the piecewise Bernoulli measure considered in \cite{Barany2018ShrinkingTO}. 
For convenience, we write
\begin{align}\label{E21}
    F(\mathtt{x},k)= \mu\big[(\mathtt{x}^{(1)}_{1:k},\mathtt{x}^{(2)}_{1: \lceil k \log_{m_2} m_1\rceil})\big],
\end{align}
and for the rest of this section abuse the notation $\mathtt{x}_i$ (similarly for $\mathtt{x}^{(1)}_i$ and $\mathtt{x}^{(2)}_i$) to denote the random variable that indicates the $i$-th letter of the elements in the sample space $(\Sigma_{m_1} \times \Sigma_{m_2})^{\mathbb{N}}$, which consequently renders $F(\mathtt{x},k)$ a random variable in this context. Under the circumstances, we recover in Lemma \ref{lem:local_dimension} an analog of \cite[Lemma 4.2]{Barany2018ShrinkingTO}:
\begin{align*}
    \liminf_{k \to \infty} \frac{-\log_{m_1} F(\mathtt{x},k)}{k} = \liminf_{k \to \infty} \mathbb{E}_{p} \left[\frac{-\log_{m_1} F(\mathtt{x},k)}{k}\right] \quad \mu\text{-a.e.},
\end{align*}
where $\mathbb{E}_{p}$ denotes the mathematical expectation with respect to $\mu$. It is noteworthy that the above limit coincides with the lower local dimension function almost everywhere (see for example \cite{K_enm_ki_2013}), i.e.,
\[
    \liminf_{k \to \infty} \frac{-\log_{m_1} F(\mathtt{x},k)}{k} =\liminf_{r\rightarrow 0} \frac{\log \mu \circ \pi^{-1}(B(\pi(\mathtt{x}),r))}{\log r} \quad \mu\text{-a.e.}
\]
We then apply \cite[Proposition 2.3]{falconer1997techniques} to conclude that 
\[
    \dim_{\mathrm{H}}(W_{\gamma}(K,T,\psi)) \ge \sup_{p} \left[ \liminf_{k \to \infty} \mathbb{E}_{p} \left[\frac{-\log_{m_1} F(\mathtt{x},k)}{k}\right] \right],
\]
for which the lower limit on the right-hand side is calculated in Lemma \ref{lem:local_dimension} with the optimal value obtained thereafter. For convenience, we introduce the following notations.
\begin{equation} \label{E22}
    H(p)=\sum\limits_{\mathtt{a} \in A} - p_{\mathtt{a}} \log_{m_1} p_{\mathtt{a}},
\end{equation}
\begin{equation} \label{E9}
    H_1(p)=\sum\limits_{\mathtt{a}^{(1)} \in A^{(1)}}- p_{\mathtt{a}^{(1)}} \log_{m_1} p_{\mathtt{a}^{(1)}},
\end{equation}
and
\begin{equation} \label{E18}
    H_2(p)=\sum\limits_{\mathtt{a} \in A}- p_{\mathtt{a}}\log_{m_1} \frac{p_{\mathtt{a}}}{p_{\mathtt{a}^{(1)}}}.
\end{equation}
 
Now, we begin to prove the lemmas.
\begin{lemma} \label{lem:Markov}
    Let $\{X_i: i \in \mathbb{N}\}$ be independent random variables with zero means and $M:=\sup_i \mathbb{E}_p |X_i|^4 < \infty$. Then, for all $n, m \in \mathbb{N}$ and $\epsilon > 0$, $\mu(|m^{-1} \sum_{i=1}^n X_i| > \epsilon) \le 3 (\frac{m^4}{n^2})^{-1} M \epsilon^{-4}$.
\end{lemma}
\begin{proof}
    It is an immediate consequence of Markov's inequality.
\end{proof}
\begin{lemma} \label{lem:local_dimension}
    Let $F(\mathtt{x},k)$ be as defined in (\ref{E21}) and suppose $0 \le \tau_1 < \infty$. Then, 
    \begin{equation}\label{lemma3.1eq}
        \liminf_{k \to \infty} \frac{-\log_{m_1} F(\mathtt{x},k)}{k} = \liminf_{k \to \infty} \mathbb{E}_{p}\left[\frac{-\log_{m_1} F(\mathtt{x},k)}{k}\right] \quad \mu\text{-a.e.}
    \end{equation}
\end{lemma}
\begin{proof}
    Recall that according to equation \eqref{eq:decomposition}, one can express $\mathtt{x}$ in the following form:
    \[
        \mathtt{x} = \mathttul{w}{1} \mathttul{u}{1} \mathttul{v}{1} \mathttul{w}{2} \mathttul{u}{2} \mathttul{v}{2} \mathttul{w}{3} \cdots.
    \]
    We first prove the claim that for any $\epsilon > 0$,
    \[
        \mu\left(\left|\frac{\log_{m_1} F(\mathtt{x},k)}{k} - \mathbb{E}_p\left(\frac{\log_{m_1} F(\mathtt{x},k)}{k}\right)\right| > \epsilon\right) = O(k^{-2}),
    \]
    where the $O(k^{-2})$ may depend on $p$ and $\epsilon$ but not on $k$, and thus $\epsilon$ throughout the discussion is treated as a constant. Due to the assumption $n_i \gg R_{i-1}$, every $k \in \mathbb{N}$ must fall in one of the following cases.
    
    \textbf{Case 1: $n_i < k$ and $\lfloor k \log_{m_2} m_1\rfloor \le L_i$.} For such $k$, we have
    \begin{equation} \label{eq:probability-1}
        \begin{aligned}
            \frac{\log_{m_1} F(\mathtt{x},k)}{k} =& \frac{1}{k}\left[\log_{m_1} \mu[\mathttul{w}{1} \cdots \mathttul{v}{i-1}] + \sum_{\ell=1}^{|\mathttul{w}{i}|} \log_{m_1} p_{\mathttul[(1)][\ell]{w}{i}} + \sum_{\ell=1}^{\min\{|\mathttul{w}{i}|, \lfloor k \log_{m_2} m_1 \rfloor - R_{i-1}\}} \log_{m_1} \frac{p_{\mathttul[][\ell]{w}{i}}}{p_{\mathttul[(1)][\ell]{w}{i}}} \right. \\
            & \qquad\left. + \sum_{\ell=1}^{\max\{k - R_i,0\}} \log_{m_1} p_{\mathttul[(1)][\ell]{w}{i+1}} \right].
        \end{aligned}
    \end{equation}
    Observe that, by virtue of $n_i \gg R_{i-1}$, 
    \begin{equation} \label{eq:o(1)}
        \left|\frac{1}{k} \log_{m_1} \mu[\mathttul{w}{1} \cdots \mathttul{v}{i-1}]\right| \le \frac{R_{i-1}}{k} \cdot \max_{\mathtt{a} \in A: p_{\mathtt{a}} \ne 0} \left|\log_{m_1} p_{\mathtt{a}}\right| = o(1).
    \end{equation}
    In addition, according to \eqref{eq:prob_measure} the remaining terms are, respectively, summations over i.i.d.~random variables
    \[
        \left\{\log_{m_1} p_{\mathttul[(1)][\ell]{w}{i}}: 1 \le \ell \le |\mathttul{w}{i}|\right\} \text{ and } \left\{\log_{m_1} \frac{p_{\mathttul[][\ell]{w}{i}}}{p_{\mathttul[(1)][\ell]{w}{i}}}: 1 \le \ell \le |\mathttul{w}{i}|\right\}
    \]
    with finite fourth moments
    \[
        \mathbb{E}_p \left|\log_{m_1} p_{\mathttul[(1)][\ell]{w}{i}}\right|^4 \le \max_{\mathtt{a}^{(1)} \in A^{(1)}: p_{\mathtt{a}^{(1)}} \ne 0} \left|\log_{m_1} p_{\mathtt{a}}\right|^4 \quad \text{and} \quad \mathbb{E}_p \left|\log_{m_1} \frac{p_{\mathttul[][\ell]{w}{i}}}{p_{\mathttul[(1)][\ell]{w}{i}}}\right|^4 \le \max_{\mathtt{a} \in A: p_{\mathtt{a}} \ne 0} \left|\log \frac{p_{\mathtt{a}}}{p_{\mathtt{a}^{(1)}}}\right|^4.
    \]
    Therefore, Lemma \ref{lem:Markov}, together with \eqref{eq:o(1)}, yields
    \[
        \mu\left(\left|\frac{\log_{m_1} F(\mathtt{x},k)}{k} - \mathbb{E}_p\left(\frac{\log_{m_1} F(\mathtt{x},k)}{k}\right)\right| > \epsilon\right) = O(k^{-2})
    \]
    
    \textbf{Case 2: $k \le n_{i+1}$ and $L_i < \lfloor k \log_{m_2} m_1\rfloor$.} We begin with a similar expression as above:
    \begin{equation} \label{eq:probability-2}
        \begin{aligned}
            \frac{\log_{m_1} F(\mathtt{x},k)}{k} =& \frac{1}{k} \left[ \log_{m_1} \mu[\mathttul{w}{1} \cdots \mathttul{v}{i-1}] + \sum_{\ell=1}^{|\mathttul{w}{i}|} \log_{m_1} p_{\mathttul[][\ell]{w}{i}} + \sum_{\ell=1}^{\min\{|\mathttul{v}{i}|, \lfloor k \log_{m_2} m_1 \rfloor - L_i\}} \log_{m_1} \frac{p_{\mathttul[][\ell]{v}{i}}}{p_{\mathttul[(1)][\ell]{v}{i}}} \right. \\
            & \qquad \left. + \sum_{\ell=1}^{\max\{k - R_i, 0\}} \log_{m_1} p_{\mathttul[(1)][\ell]{w}{i+1}} + \sum_{\ell=1}^{\max\{k \log_{m_2} m_1 - R_i, 0\}} \frac{p_{\mathttul[][\ell]{w}{i+1}}}{p_{\mathttul[(1)][\ell]{w}{i+1}}} \right].
        \end{aligned}
    \end{equation}
    The argument for i.i.d.~random variables proceeds as above except for the summation involving $\mathttul{v}{i}$, which are not necessary i.i.d. However, this issue could be resolved by considering the conditional measure $\mu_{\mathtt{w}^*}$ of $\mu$ on the cylinder set $[\mathtt{w}^*]$, $\mathtt{w}^* = \mathttul{w}{1} \mathttul{u}{1} \cdots \mathttul{w}{i} \in \mathcal{W}_i$. Hence, with respect to any well-defined $\mu_{\mathtt{w}^*}$,
    \[
        \sum_{\ell=1}^{n} \log_{m_1} \frac{p_{\mathttul[][\ell]{v}{i}}}{p_{\mathttul[(1)][\ell]{v}{i}}}
    \]
    is a summation over independent random variables with uniformly bounded fourth moments
    \[
        \mathbb{E}_p \left|\log_{m_1} \frac{p_{\mathttul[][\ell]{v}{i}}}{p_{\mathttul[(1)][\ell]{v}{i}}}\right|^4 \le \max_{\mathtt{a} \in A: p_{\mathtt{a}} \ne 0} \left|\log_{m_1} \frac{p_{\mathtt{a}}}{p_{\mathtt{a}^{(1)}}}\right|^4
    \]
    Applying Lemma \ref{lem:Markov} again yields that for every $\mathtt{w}^* \in \mathcal{W}_i$,
    \[
        \mu_{\mathtt{w}^*}\left(\left|\frac{1}{k} \sum_{\ell=1}^{n} \log_{m_1} \frac{p_{\mathttul[][\ell]{v}{i}}}{p_{\mathttul[(1)][\ell]{v}{i}}}-\mathbb{E}_p\left(\frac{1}{k} \sum_{\ell=1}^{n} \log_{m_1} \frac{p_{\mathttul[][\ell]{v}{i}}}{p_{\mathttul[(1)][\ell]{v}{i}}} \middle| \mathtt{w}^*\right)\right| > \epsilon\right) = O(k^{-2}).
    \]
    We note that the bound $O(k^{-2})$ could be chosen uniform over $\mathtt{w}^*$, and thus it remains to show that
    \[
        \mu\left(\left|\mathbb{E}_p\left(\frac{1}{k} \sum_{\ell=1}^{n} \log_{m_1} \frac{p_{\mathttul[][\ell]{v}{i}}}{p_{\mathttul[(1)][\ell]{v}{i}}} \right) - \mathbb{E}_p\left(\frac{1}{k} \sum_{\ell=1}^{n} \log_{m_1} \frac{p_{\mathttul[][\ell]{v}{i}}}{p_{\mathttul[(1)][\ell]{v}{i}}} \middle| \mathtt{w}^*\right)\right| > \epsilon\right) = O(k^{-2}).
    \]
    Recall that equation \eqref{eq:prob_measure} asserts $\mu$-almost surely,
    \[
        (\mathttul{u}{i} \mathttul{v}{i})^{(1)} = \left(\mathttul[(1)]{w}{1} \cdots \mathttul[(1)]{u}{i-1} \mathttul[(1)]{v}{i-1} \mathttul[(1)]{w}{i}\right)^{\frac{\hat{\ell}_{1}(n_i)}{n_i}},
    \]
    which implies that if we associate each index $\ell \in [1, |\mathttul{v}{i}|]$ of $\mathttul{v}{i}$ with $\ell' \in [1, n_i]$ satisfying $\ell' = L_{i} + \ell \pmod{n_i}$, then 
    \begin{equation*}
        \mathttul[][\ell]{v}{i} = \mathttul[][\ell'-R_{i-1}]{w}{i} \text{ whenever } R_{i-1} < \ell' \le n_i.
    \end{equation*}
    Due to $n_i \gg R_{i-1}$, we again deduce that for every $1 \le n \le |\mathttul{v}{i}|$,
    \[
        \frac{1}{k} \sum_{\ell=1}^{n} \log_{m_1} \frac{p_{\mathttul[][\ell]{v}{i}}}{p_{\mathttul[(1)][\ell]{v}{i}}} = \frac{1}{k} \sum_{\substack{\ell \in [1,|\mathttul{v}{i}|]: \\ \ell' \in (R_{i-1},n_i]}} \log_{m_1} \frac{p_{\mathttul[][\ell]{v}{i}}}{p_{\mathttul[(1)][\ell'-R_{i-1}]{w}{i}}} + o(1),
    \]
    where $o(1)$ is uniform due to a similar argument as \eqref{eq:o(1)}. Thus, writing
    $$
        t_{n, \ell} := \# \{j \in [1,n]: \ell = j \ (\bmod\ {n_i})\} \in \left\{\left\lfloor\frac{n}{n_i}\right\rfloor, \left\lfloor\frac{n}{n_i}\right\rfloor+1\right\} \qquad (1 \le \ell \le |\mathtt{w}^{[i]}|),
    $$
    we have
    \begin{align} \label{eq:conditional_expectation_expression}
        & \mathbb{E}_p\left(\frac{1}{k} \sum_{\ell=1}^{n} \log_{m_1} \frac{p_{\mathttul[][\ell]{v}{i}}}{p_{\mathttul[(1)][\ell]{v}{i}}} \middle| \mathtt{w}^*\right) = \frac{1}{k} \sum_{\ell=1}^{|\mathttul{w}{i}|} t_{n, \ell} \sum_{\substack{\mathtt{a} \in A: \\ \mathtt{a}^{(1)} = \mathttul[(1)][\ell]{w}{i}}}\frac{p_{\mathtt{a}_{\ell}}}{p_{\mathttul[(1)][\ell]{w}{i}}} \log_{m_1} \frac{p_{\mathtt{a}_{\ell}}}{p_{\mathttul[(1)][\ell]{w}{i}}} + o(1),
    \end{align}
    which again is a summation of i.i.d.~random variables and Lemma \ref{lem:Markov} applies. In summary, we have
    \[
        \mu\left(\left|\frac{\log_{m_1} F(\mathtt{x},k)}{k} - \mathbb{E}_p\left(\frac{\log_{m_1} F(\mathtt{x},k)}{k}\right)\right| > \epsilon\right) = O(k^{-2}).
    \]
    
    To wrap up, we note the lemma follows from the Borel-Cantelli lemma applied to our claim with $\epsilon \to 0$.
\end{proof}

\begin{lemma}\label{L5}
    Let $ F(\mathtt{x},k)$ be defined in (\ref{E21}). Then $\liminf_{k \to \infty} \mathbb{E}_{p} \left[\frac{-\log_{m_1} F(\mathtt{x},k)}{k}\right] = \min\{t_1, t_2\}$, where
    \[
        t_1 = \begin{cases}
            \frac{(1 - (1 - \gamma) \tau_1 \log_{m_2} m_1) \cdot H_1(p) + \log_{m_2} m_1 \cdot H_2(p)}{1+\tau_2} & \text{if } 1+(1-\gamma)\tau_1 \le \log_{m_1} m_2, \\
            \frac{\log_{m_2} m_1 \cdot H_1(p) + \log_{m_2} m_1 \cdot H_2(p)}{1+\tau_2} & \text{if } 1+(1-\gamma)\tau_1 > \log_{m_1} m_2.
        \end{cases}
    \]
    \[
        t_2 = \begin{cases}
            \frac{H_1(p)}{1+\tau_1} + \log_{m_2} m_1 \cdot H_2(p) & \text{if } 1 + \tau_1 \le \log_{m_1} m_2, \\
            \frac{H_1(p) + H_2(p)}{1+\tau_1} & \text{if } 1 + (1 - \gamma) \tau_1 \le \log_{m_1} m_2 < 1 + \tau_1, \\
            \frac{H_1(p) + \log_{m_2} m_1 \cdot (1 + (1 - \gamma) \tau_1) \cdot H_2(p)}{1+\tau_1} & \text{if } 1+(1-\gamma) \tau_1 > \log_{m_1} m_2.
        \end{cases}
    \]
\end{lemma}
\begin{proof}
    Recall that $H(p) = H_1(p) + H_2(p)$. Calculating expectation of \eqref{eq:probability-1} and \eqref{eq:probability-2} respectively yields that
    
    \textbf{Case 1:} if $n_i < k \le L_i \log_{m_1} m_2$, then
    \begin{equation*}
        \begin{aligned}
            \mathbb{E}_p\left[\frac{\log_{m_1} F(\mathtt{x},k)}{k}\right] &= \frac{1}{k}\left[(n_i-R_{i-1}+\max\{k - R_i,0\}) H_1(p) \right. \\
            & \left. \qquad + (n_i-R_{i-1}+\min\{0, k \log_{m_2} m_1 - n_i\}) H_2(p) \right] + o(1);
        \end{aligned}
    \end{equation*}
    
    \textbf{Case 2:} if $L_i \log_{m_1} m_2 < k \le n_i$, then
    \begin{equation*}
        \begin{aligned}
            \mathbb{E}_p\left[\frac{\log_{m_1} F(\mathtt{x},k)}{k}\right] &= \frac{1}{k} \left[ (n_i-R_{i-1}+\max\{k - R_i, 0\}) H_1(p) \right.\\
            & \left. \qquad + (n_i-R_{i-1} + k \log_{m_2} m_1 - L_i) H_2(p) \right] + o(1).
        \end{aligned}
    \end{equation*}
    From these expressions, we observe that the interval $(n_i, n_{i+1}]$ can be divided by $\{L_i, R_i, n_i \log_{m_1} m_2, L_i \log_{m_1} m_2\}$ into at most five subintervals so that on each subinterval the function $k \mapsto \mathbb{E}_p\left[\frac{\log_{m_1} F(\mathtt{x},k)}{k}\right]$, by neglecting the $o(1)$ terms, is a polynomial of $k$ of the form $a k^{-1}+b$ and hence monotone. This simplifies the calculation of the desired lower limit down to evaluating the minimum over the numbers $\{S_1, S_2, S_3, S_4, S_5\}$, which are defined to be $\liminf_{i \to \infty} \mathbb{E}_p\left[\frac{\log_{m_1} F(\mathtt{x},k_i)}{k_i}\right]$ with $k_i$ taken as $n_i$, $L_i$, $R_i$, $n_i \log_{m_1} m_2$, and $L_i \log_{m_1} m_2$, respectively. These values could be straightforwardly calculated and are summarized in Table \ref{tab:S_i}.
    \begin{table}[]
        \centering
        \begin{tabular}{ccc}
            $k_i$ & Case & $\liminf_{i \to \infty} \mathbb{E}_p\left[\frac{\log_{m_1} F(\mathtt{x},k_i)}{k_i}\right]$ \\
            \hline
            \multirow{ 2}{*}{$n_i$} & 1 & N/A \\
            & 2 & $H_1(p) + \log_{m_2} m_1 \cdot H_2(p)$ \\
            \hline
            \multirow{ 2}{*}{$L_i$} & 1 & $(1+\tau_2)^{-1} \cdot H_1(p) + \max\{(1+\tau_2)^{-1},\log_{m_2} m_1\} \cdot H_2(p)$ \\
            & 2 & N/A \\
            \hline
            \multirow{ 2}{*}{$n_i \log_{m_1} m_2$} & 1 & $\max\{1 - \tau_1 \log_{m_2} m_1, \log_{m_2} m_1\} \cdot H_1(p) + \log_{m_2} m_1 \cdot H_2(p)$ \\
            & 2 & N/A \\
            \hline
            \multirow{ 2}{*}{$L_i \log_{m_1} m_2$} & 1 & $\max \{1 - (1+\tau_2)^{-1} \tau_1 , (1+\tau_2)^{-1} \} \cdot H_1(p)\log_{m_2} m_1 + (1+\tau_2)^{-1} \log_{m_2} m_1 \cdot H_2(p)$ \\
            & 2 & N/A \\
            \hline
            \multirow{ 2}{*}{$R_i$} & 1 & $(1+\tau_1)^{-1} \cdot H_1(p) + \min\{(1+\tau_1)^{-1},\log_{m_2} m_1\} \cdot H_2(p)$ \\
            & 2 & $(1+\tau_1)^{-1} \cdot H_1(p) + (\log_{m_2} m_1 - (1+\tau_1)^{-1} \tau_2) \cdot H_2(p)$
        \end{tabular}
        \caption{values of $S_i$}
        \label{tab:S_i}
    \end{table}
    Since $H_1(p), H_2(p), \tau_1, \tau_2 \ge 0$, we observe the following order among these values.
    \begin{itemize}
        \item $S_1 \ge S_3 \ge S_5$. The inequality $S_1 \ge S_3$ is clear by comparing the coefficients of $H_1(p)$ and of $H_2(p)$. As for $S_3 \ge S_5$, it is done by comparing the coefficients as follows.
        \[
            \max\{1 - \tau_1 \log_{m_2} m_1, \log_{m_2} m_1\} = \begin{cases}
                \log_{m_2} m_1 \ge (1+\tau_1)^{-1} & \text{if } \log_{m_1} m_2 \le 1 + \tau_1, \\
                1 - \tau_1 \log_{m_2} m_1 < (1+\tau_1)^{-1} & \text{if } \log_{m_1} m_2 \ge 1 + \tau_1,
        \end{cases}
        \]
        \[
            \log_{m_2} m_1 \ge \min\{(1+\tau_1)^{-1}, \log_{m_2} m_1\} \quad{and}\quad \log_{m_2} m_1 - (1+\tau_1)^{-1}\tau_2.
        \]
        \item $S_2 \ge S_4$ if $\tau_1 \ge \tau_2$. Note that the coefficient of $H_1(p)$ in $S_4$ could be viewed as the maximum of two linear functions of $\log_{m_2} m_1$, which yields
        \[
            \max \{1 - (1+\tau_2)^{-1} \tau_1 \log_{m_2} m_1, (1+\tau_2)^{-1} \log_{m_2} m_1\} \le (1+\tau_1)^{-1}.
        \]
        Under our assumption, this is less than the corresponding coefficient $(1+\tau_2)^{-1}$ in $S_2$. The order between the coefficients of $H_2(p)$ is clear.
    \end{itemize}
    Now we plug $\tau_2 = \gamma \log_{m_2} m_1 \cdot \tau_1$ into $S_4$ and $S_5$ to deduce
    \[
        S_4 = \begin{cases}
            \frac{(1 - (1 - \gamma) \tau_1 \log_{m_2} m_1) \cdot H_1(p) + \log_{m_2} m_1 \cdot H_2(p)}{1+\tau_2} & \text{if } 1+(1-\gamma)\tau_1 \le \log_{m_1} m_2, \\
            \frac{\log_{m_2} m_1 \cdot H_1(p) + \log_{m_2} m_1 \cdot H_2(p)}{1+\tau_2} & \text{if } 1+(1-\gamma)\tau_1 > \log_{m_1} m_2.
        \end{cases}
    \]
    \[
        S_5 = \begin{cases}
            \frac{H_1(p)}{1+\tau_1} + \log_{m_2} m_1 \cdot H_2(p) & \text{if } 1 + \tau_1 \le \log_{m_1} m_2, \\
            \frac{H_1(p) + H_2(p)}{1+\tau_1} & \text{if } 1 + (1 - \gamma) \tau_1 \le \log_{m_1} m_2 < 1 + \tau_1, \\
            \frac{H_1(p) + \log_{m_2} m_1 \cdot (1 + (1 - \gamma) \tau_1) \cdot H_2(p)}{1+\tau_1} & \text{if } 1+(1-\gamma) \tau_1 > \log_{m_1} m_2.
        \end{cases}
    \]
    The proof is then completed.
\end{proof}
\begin{remark} \label{rem:comparison_barany}
    The essential difference between the measure $\mu$ constructed in this section and its counterpart in \cite{Barany2018ShrinkingTO} lies in the subwords $\mathttul[]{u}{j}$ and $\mathttul[]{v}{j}$ ($1 \le j \le i$). Specifically, these subwords satisfy \eqref{eq:recurrence_assumption} in our work, while they fulfill a priori assumptions imposed on the shrinking target problem in \cite{Barany2018ShrinkingTO}. The nuance urges necessary modifications in the proof of constant lower limit (Lemma \ref{lem:local_dimension}) and is reflected in the expression \eqref{eq:conditional_expectation_expression}. However, the subtlety is subsumed in $\mathbb{E}_{p} \left[\frac{-\log_{m_1} F(\mathtt{x},k)}{k}\right]$ and in the homogeneity assumption of the Bedford-McMullen carpet; therefore, it is barely visible beyond the scope of the proof of Lemma \ref{lem:local_dimension}.
\end{remark}
Finally, we wrap up this section by noting that the lower bound of $\dim_{\mathrm{H}} W_{\gamma}(K,T,\psi)$ in Theorem \ref{M1-1} holds as a consequence of Lemma \ref{lem:local_dimension} and \ref{L5}, for which we simply take $p_{\mathtt{a}}=\frac{1}{\# A}$ for all $\mathtt{a} \in A$ so that $H_1(p)=\log_{m_1} M$ and $H_2(p)=\log_{m_1} N$.

\section{Proof of the upper bound}\label{S3}
    In this section, we give the rest of the proof of Theorem \ref{M1}. We note that it suffices to prove the case where $\lim_{n \to \infty}-\frac{\log \psi(n)}{n}$ exists, since for those $\psi(n)$ without this property, we may consider $\phi(n):=\max\{\psi(n),m_1^{-\tau_1 \cdot n}\}$ so that $\dim_{\mathrm{H}} W_{\gamma}(K,T,\psi) \le \dim_{\mathrm{H}} W_{\gamma}(K,T,\phi)$ and we may apply the aforementioned result of the case when $\lim_{n \to \infty}-\frac{\log \psi(n)}{n}$ exists.
    
    Note that the set $W_{\gamma}(K,T,\psi)$ can be written as a limsup set 
    \[
        W_{\gamma}(K,T,\psi)=\limsup_{n \to \infty} W_n^{\gamma}(K,T,\psi),
    \]
    where
    \[
        W_n^{\gamma}(K,T,\psi):=\left\{x \in K: \begin{cases}
            |x^{(1)}-T^n_{m_1}(x^{(1)})| < \psi(n) \\
            |x^{(2)}-T^n_{m_2}(x^{(2)})| < \psi(n)^{\gamma} 
        \end{cases}\right\}.
    \]
    We note that the set above can further be written as the following union:
    \[
        W_n^{\gamma}(K,T,\psi)=\bigcup_{\mathtt{w} \in A^n} J_n(\mathtt{w}) := \bigcup_{\mathtt{w} \in A^n} W_{n}^{\gamma}(K,T,\psi) \cap I(\mathtt{w}).
    \]
    For sufficiently large $n$, the set $J_n(\mathtt{w})$ is contained in the interior of a rectangle $R^{(1)}(\mathtt{w})\times R^{(2)}(\mathtt{w})$ of width $4\psi(n)m_1^{-n}$ and height $4\psi(n)^{\gamma}m_2^{-n}$, since for any $\pi(\mathtt{x}) \in J_n(\mathtt{w})$,
    \begin{equation} \label{eq:recurrence_set_diameter}
        \begin{aligned}
            |\pi(\mathtt{x})^{(i)}-\pi(\mathtt{w}^{\infty})^{(i)}| &\ge |T^n(\pi(\mathtt{x}))^{(i)}-T^n(\pi(\mathtt{w}^{\infty}))^{(i)}|-|\pi(\mathtt{x})^{(i)}-T^n(\pi(\mathtt{x})^{(i)})| \\
            &\ge m_i^n |\pi(\mathtt{x})^{(i)}-\pi(\mathtt{w}^{\infty})^{(i)}| - |\pi(\mathtt{x})^{(i)}-T^n(\pi(\mathtt{x})^{(i)})|,
        \end{aligned}
    \end{equation}
    $|\pi(\mathtt{x})^{(1)}-T^n(\pi(\mathtt{x})^{(1)})| \leq \psi(n)m_1^{-n}$ and $|\pi(\mathtt{x})^{(1)}-T^n(\pi(\mathtt{x})^{(1)})|\leq \psi(n)^{\gamma}m_2^{-n}$.
    By Remark \ref{remark-1}, we know that $0<\tau_1<\infty$, then there exists $G_1\in \mathbb{N}$ such that for any $n\geq G$, $\psi(n)<1, (\frac{m_2}{m_1})^{n}>4.$
    Now for any given $\delta>0$, we can choose $G$ large enough so that for any $n\geq \max\{G,G_1\}$, $\max\{4\psi(n)m_i^{-n}, 4\psi(n)^{\gamma}m_2^{-n}\}<\delta$.
    For any $n\geq 1$, let
    \begin{equation}
        \psi_1(n)=\psi(n) \text {\ and \ }\psi_2(n)=\psi(n)^{\gamma}.
    \end{equation}
    For $i=1,2$,
    we denote by $\mathcal{B}_{i,n}$ the smallest number of balls of diameter $4\psi_i(n)m_i^{-n}$
    (the sidelength of the rectangles in $R^{(1)}(\mathtt{w})\times R^{(2)}(\mathtt{w})$ along the direction of the $i$-th axis) needed to cover $W_n^{\gamma}(K,T,\psi)$, the $s$-dimensional Hausdorff measure has the following estimate:
    \begin{align*}
        \mathcal{H}_{\delta}^s(W_{\gamma}(K,T,\psi))\leq \sum_{n=G}^{\infty} \# \mathcal{B}_{i,n}  \cdot (4\psi_i(n)m_i^{-n})^s.
    \end{align*}

We now estimate the value of $\# \mathcal{B}_{i,n}$, $i=1,2$.

For $i=2$, we need to estimate the number of balls $\mathcal{B}_{2,n}$ of diameter $4\psi(n)^{\gamma}m_2^{-n}$ to cover $W_n^{\gamma}(K,T,\psi)=\cup_{\mathtt{w} \in A^n} J_n(\mathtt{w})$. Since $J_n(\mathtt{w})$ is contained in the interior of a rectangle $R^{(1)}(\mathtt{w})\times R^{(2)}(\mathtt{w})$ of width $4\psi(n)m_1^{-n}$ and height $4\psi(n)^{\gamma}m_2^{-n}$. To estimate the values, we consider two possible cases depending on the values of $(1-\gamma)\tau_1$ and $\log_{m_1}m_2-1$. 

\textbf{Case a-1:} $(1-\gamma)\tau_1>\log_{m_1}m_2-1$. By the definition of $\tau_1$, $4\psi(n)^{\gamma}m_2^{-n}>4\psi(n)m_1^{-n}$. For any $n\geq G_1$, $m_1^{-n}\geq 4\psi(n)^{\gamma}m_2^{-n}$ by the fact that $\psi(n)< 1$ and $m_1\leq m_2$.
It is clear that $\#\mathcal{B}_{2,n} \le (MN)^n$. Therefore,
\begin{align}\label{E15}
    \mathcal{H}_{\delta}^s(W_{\gamma}(K,T,\psi))&\leq \sum_{n=G}^{\infty} (MN)^n (4\psi(n)^{\gamma}m_2^{-n})^s\\\nonumber
    &= C\sum_{n=G}^{\infty}m_1^{n h_n},
\end{align}
where $C:=4^s$ is a constant and
\begin{align*}
    h_n=\log_{m_1}M+\log_{m_1}N+s(n^{-1}\gamma\log_{m_1}\psi(n)-\log_{m_1}m_2).
\end{align*}
Since we assume $\lim_{n \to \infty}-\frac{\log \psi(n)}{n}$ exists, it is clear that $\sum_{n=G}^{\infty}m_1^{n \cdot h_n}$ converges as long as $\lim_{n \to \infty}h_n<0$ and that this is equivalent to the condition that 
\begin{align*}
    s> \lim_{n \to \infty}\frac{\log_{m_1}M+\log_{m_1}N}{\log_{m_1}m_2-n^{-1}\gamma\log_{m_1}\psi(n)} = \frac{\log_{m_2}M+\log_{m_2}N}{1+\tau_2}.
\end{align*}
This together with \eqref{E15} we can get that 
\begin{align*}
    0\leq \mathcal{H}^s(W_{\gamma}(K,T,\psi))=\lim_{\delta \to 0} \mathcal{H}_{\delta}^s(W_{\gamma}(K,T,\psi))\leq \lim_{G \to \infty}C \sum_{n=G}^{\infty} m_1^{n \cdot h_n}=0,
\end{align*}
which implies that $\dim_{\mathrm{H}}W_{\gamma}(K,T,\psi)\leq \frac{\log_{m_2}M+\log_{m_2}N}{1+\tau_2}$.

\textbf{Case a-2:} $(1-\gamma)\tau_1\leq \log_{m_1}m_2-1$. By the definition of $\tau_1$, $4\psi(n)^{\gamma}m_2^{-n}\leq 4\psi(n)m_1^{-n}$.
Notice that just a ball of diameter $4\psi(n)^{\gamma}m_2^{-n}$ is needed to cover $J_n(\mathtt{w})$ along the direction of the second axis.
Now let us look in the direction along the first axis.
Let $u>0$ be the unique integer such that $m_1^{-u}\leq 4\psi(n)m_1^{-n} \leq m_1^{-u+1}$, which implies there exist $\mathtt{w}'$, $\mathtt{w}''\in A^{u-1}$ such
that
\begin{align}\label{E16}
    J_n(\mathtt{w})\subseteq (I^{(1)}_{u-1}(\mathtt{w}')\times R^{(2)}(\mathtt{w}) \cup I^{(1)}_{u-1}(\mathtt{w}'')\times R^{(2)}(\mathtt{w}))\cap K,
\end{align}
where
$I^{(1)}_{u-1}(\mathtt{w}')$ is the projection of the rectangle $I_{u-1}(\mathtt{w}')$
onto the first axis, 
$I^{(1)}_{u-1}(\mathtt{w}'')$ is that of $I_{u-1}(\mathtt{w}'')$, and 
$K$ is the Bedford-McMullen carpet.
Since $4\psi(n)^{\gamma}m_2^{-n}\leq 4\psi(n)m_1^{-n}$, there exists a unique integer $v\geq 0$ such that
\begin{align}\label{E4}
    m_1^{-u-v}\leq 4\psi(n)^{\gamma}m_2^{-n} \leq m_1^{-u-v+1},
\end{align}
which yields that $v\leq 1-n+n\log_{m_1}m_2+(1-\gamma)\log_{m_1}\psi(n)$.
Then $J_n(\mathtt{w})$ can be covered by $2M^{v+1}$ balls of diameter $4\psi(n)m_2^{-n}$, which 
follows from that $M$ denotes the number of columns containing at least one chosen rectangle. By  \eqref{E16}, \eqref{E4} and the definition of $ \mathcal{B}_{2,n}$,
\[ 
    \# \mathcal{B}_{2,n}\leq 2 M^{v+1} \cdot (\#A)^n= 2M^{v+1} (MN)^n.
\]
The remaining argument is just the same as \textbf{Case a-1}. We therefore conclude that
\[
    \dim_{\mathrm{H}}W_{\gamma}(K,T,\psi)\leq (1-\frac{\tau_1\log_{m_2}m_1}{1+\tau_2})\log_{m_1}M+\frac{\log_{m_2}N}{1+\tau_2}.
\]

For $i=1$, we need to estimate the number of balls $\mathcal{B}_{1,n}$ of diameter $4\psi(n)m_1^{-n}$ to cover $W_n(K,T,\psi)=\cup_{\mathtt{w} \in A^n} J_n(\mathtt{w})$. To estimate the values, we consider two possible cases depending on the values of $\log_{m_1} m_2 $ and $1+\tau_1$.

\textbf{Case b-1:} $\log_{m_1} m_2 >1+\tau_1$. By the definition of $\tau_1$, $4\psi(n)m_1^{-n}>m_2^{-n}$. Now observe that a ball of diameter $4\psi(n)m_1^{-n}$ will also cover other $n$-th level rectangles in $\bigcup_{\mathtt{w} \in A^n} R^{(1)}(\mathtt{w})\times R^{(2)}(\mathtt{w})$ along the second axis. Let $j=\lfloor \log_{m_2}\psi(n)-n\log_{m_2}m_1+n\rfloor-2$, then $\#\mathcal{B}_{1,n}\leq (MN)^n N^{-j}$, which yields that
\begin{align*}
    \mathcal{H}_{\delta}^s(W_{\gamma}(K,T,\psi))\leq  \sum_{n=G}^{\infty} (MN)^n N^{-j} (4\psi(n)m_1^{-n})^s.
\end{align*}
Such an argument also applies to the case $i=1$, from which we conclude that $\dim_{\mathrm{H}}W_{\gamma}(K,T,\psi)\leq \log_{m_2}N+\frac{\log_{m_1}M}{1+\tau_1}$.

\textbf{Case b-2:} $\log_{m_1} m_2 \leq 1+\tau_1$. By the definition of $\tau_1$, $4\psi(n)m_1^{-n}\leq m_2^{-n}$. 
    
\textbf{Case b-2-1:} $(1-\gamma)\tau_1<\log_{m_1}m_2-1$. By the definition of $\tau_1$, $4\psi(n)^{\gamma}m_2^{-n}<4\psi(n)m_1^{-n}$.
It is clear that $\#\mathcal{B}_{1,n} \le (MN)^n$. Therefore,
\begin{align*}
    \mathcal{H}_{\delta}^s(W_{\gamma}(K,T,\psi)) \leq \sum_{n=G}^{\infty} (MN)^n (4\psi(n)m_1^{-n})^s.
\end{align*}
The remaining argument is just the same as the case when $i=1$. We therefore conclude that
\[
    \dim_{\mathrm{H}}W_{\gamma}(K,T,\psi)\leq \frac{\log_{m_1}M+\log_{m_1}N}{1+\tau_1}.
\]
\textbf{Case b-2-2:} $(1-\gamma)\tau_1> \log_{m_1}m_2-1$. By the definition of $\tau_1$, $4\psi(n)^{\gamma}m_2^{-n}\geq 4\psi(n)m_1^{-n}$. Let $u>0$ be the unique integer such that $m_1^{-u}\leq 4\psi(n)^{\gamma}m_2^{-n} \leq  m_1^{-u+1}$, which implies there exist $\mathtt{w}'$, $\mathtt{w}''\in A^{u-1}$ such that
\begin{align}\label{E16-1}
    J_n(\mathtt{w})\subseteq (I^{(1)}_{u-1}(\mathtt{w}')\times R^{(2)}(\mathtt{w}) \cup I^{(1)}_{u-1}(\mathtt{w}'')\times R^{(2)}(\mathtt{w}))\cap K,
\end{align}
where $I^{(1)}_{u-1}(\mathtt{w}')$ is the projection of the rectangle $I_{u-1}(\mathtt{w}')$ onto the first axis, $I^{(1)}_{u-1}(\mathtt{w}'')$ is that of $I_{u-1}(\mathtt{w}'')$, and $K$ is the Bedford-McMullen carpet. Since $4\psi(n)^{\gamma}m_2^{-n}\geq 4\psi(n)m_1^{-n}$, there exists a unique integer $v\geq 0$ such that
\begin{align}\label{E4-1}
    m_1^{-u-v}\leq 4\psi(n)m_1^{-n} \leq m_1^{-u-v+1},
\end{align}
which yields that $v\leq 1-n+n\log_{m_2} m_1+(\gamma-1)\log_{m_2}\psi(n)$. That $J_n(\mathtt{w})$ can be covered by $2N^{v+1}$ balls of diameter $4\psi(n)m_1^{-n}$. By \eqref{E16-1}, \eqref{E4-1} and the definition of $ \mathcal{B}_{1,n}$,
\[
    \# \mathcal{B}_{1,n}\leq 2 N^{v+1} \cdot (\#A)^n= 2N^{v+1} (MN)^n.
\]
The remaining argument is just the same as \textbf{Case a-1}. We therefore conclude that
\[
    \dim_{\mathrm{H}}W_{\gamma}(K,T,\psi)\leq \frac{(1+(1-\gamma)\tau_1)\log_{m_2}N+\log_{m_1}M}{1+\tau_1}.
\]

From the discussion above, we obtain the upper bound for Theorem \ref{M1-1}.

\section{Applications}\label{S4}
In this section, we present some examples. Example \ref{T3} is a special case of Theorem \ref{M1} for which $m_1=M$ and $m_2=N$.
\begin{example}\label{T3}
    Let $T: [0,1]^2 \to [0,1]^2 $ be an integer diagonal matrix transformation of $[0,1]^2$, i.e., $T(x)=(T_{m_1}(x^{(1)}), T_{m_2}(x^{(2)})):=\left(m_1 x^{(1)}\ (\bmod\ 1),m_2 x^{(2)}\ (\bmod\ 1)\right)$, with $2\leq m_1\leq m_2$. Suppose $\psi: \mathbb{N}^{+}\to \mathbb{R}^{+}$ be a real positive function. Let 
    \begin{align*}
        W(T,\psi):=\left\{\mathtt{x}\in [0,1]^2:
        T^n(x) \in B(x, \psi(n)) \ \text{for infinitely many } n.
        \right\}
    \end{align*}
    Then, 
    \begin{align*}
        \dim_{\mathrm{H}} W(T,\psi) = \begin{cases}
            \min \left\{\frac{2}{1+\tau_2}, \frac{1}{1+\tau_1} + 1\right\} & \text{if } \log_{m_1} m_2 > 1 + \tau_1, \\
            \min\left\{\frac{2}{1+\tau_2}, \frac{1 +\log_{m_1} m_2}{1+\tau_1} \right\} & \text{if } \log_{m_1} m_2 \le 1 + \tau_1,
        \end{cases}
    \end{align*}
    where $\tau_{i}$, $i=1,2$, is defined in (\ref{E1}).
\end{example}
The following example illustrates the case when $K$ is a product of Cantor sets.
\begin{example}
    Let $\mathcal{C}_{\frac{1}{3}} $ denote the middle third Cantor set and $\mathcal{C}_{\frac{1}{4}} $ be the attractor of the iterated function system $\{f_1,f_2,f_3\}$ on $[0,1]$, where $ f_1(x)=\frac{1}{4}x, f_2(x)=\frac{1}{4}x+\frac{1}{4}$ and $f_3(x)=\frac{1}{4}x+\frac{3}{4}$. Define $T: \mathcal{C}_{\frac{1}{3}}\times \mathcal{C}_{\frac{1}{4}} \to \mathcal{C}_{\frac{1}{3}}\times \mathcal{C}_{\frac{1}{4}}$ as 
    \[
        T(x)=\big(T_{3}(x^{(1)}), T_{4}(x^{(2)})\big):=\left(3 x^{(1)}\ (\bmod\ 1), 4 x^{(2)}\ (\bmod\ 1)\right).
    \]
    Suppose $\psi: \mathbb{N}^{+}\to \mathbb{R}^{+}$ is a real positive function.
    Let 
    \begin{align*}
        W_1(T,\psi):=\left\{\mathtt{x}\in \mathcal{C}_{\frac{1}{3}}\times \mathcal{C}_{\frac{1}{4}}: T^n(x) \in B(x, \psi(n)) \ \text{for infinitely many } n \right\}.
    \end{align*}
    Then,
    \begin{align*}
        \dim_{\mathrm{H}} W_1(T,\psi) = \begin{cases}
            \min \left\{\frac{\log_{3} 2 + \log_{4} 3}{1+\tau_2}, \frac{\log_{3} 2}{1+\tau_1} + \log_{4} 3\right\} & \text{if } \log_{3} 4 > 1 + \tau_1, \\
            \min\left\{\frac{\log_{3} 2 + \log_{4} 3}{1+\tau_2}, \frac{\log_{3} 2 +1}{1+\tau_1} \right\} & \text{if } \log_3 4 \le 1 + \tau_1,
        \end{cases}
    \end{align*}
    where $\tau_{i}$, $i=1,2$, is defined in (\ref{E1}).
\end{example}
\textbf{Acknowledgements.} The authors would like to thank Professor Bing Li for his guidance and suggestions. The authors would like to thank Professor Meng Wu for drawing our attention to these problems and reading the early draft of the manuscript. The authors thank the referees for their valuable and very detailed suggestions.

\bibliographystyle{amsplain}

\end{document}